\documentclass[MSNbibl,number,citesort,rotating,dvips]{arxbj}
\usepackage{upgreek,multirow}
\usepackage[left]{vmkcol}
\usepackage{graphicx}

\aid{0}
\volume{18}
\issue{4}
\pubyear{2012}
\firstpage{1448}
\lastpage{1464}
\doi{10.3150/11-BEJ374}

\makeatletter

\newcommand{\updown}{\bigtriangledown}
\newremark{remark}{Remark}
\newtheorem{theorem}{Theorem}

\newcolumntype{d}[1]{D{.}{.}{#1}}

\newcommand{\no}{\nonumber}
\newcommand{\om}{\omega}

\makeatother

\begin{document}
\begin{frontmatter}

\title{Inference of seasonal long-memory aggregate time series}

\runtitle{Seasonal long-memory time series}

\begin{aug}
\author[1]{\fnms{Kung-Sik} \snm{Chan}\corref{}\thanksref{1}\ead[label=e1]{kung-sik-chan@uiowa.edu}}
\and
\author[2]{\fnms{Henghsiu} \snm{Tsai}\thanksref{2}\ead[label=e2]{htsai@stat.sinica.edu.tw}}

\runauthor{K.-S. Chan and H. Tsai}
\address[1]{Department of Statistics \& Actuarial Science,
University of Iowa, Iowa City, IA 52242, USA.\\ \printead{e1}}
\address[2]{Institute of Statistical Science, Academia Sinica, Taipei
115, Taiwan, Republic of China.\\ \printead{e2}}
\end{aug}

\received{\smonth{4} \syear{2010}}
\revised{\smonth{1} \syear{2011}}

%
\begin{abstract}
Time-series data with regular and/or seasonal long-memory are often
aggregated before analysis. Often, the aggregation scale is large
enough to remove any short-memory components of the underlying process
but too short to eliminate seasonal patterns of much longer periods. In
this paper, we investigate the limiting correlation structure of
aggregate time series within an intermediate asymptotic framework that
attempts to capture the aforementioned sampling scheme. In particular,
we study the autocorrelation structure and the spectral density
function of aggregates from a discrete-time process. The underlying
discrete-time process is assumed to be a stationary Seasonal
AutoRegressive Fractionally Integrated Moving-Average (SARFIMA)
process, after suitable number of differencing if necessary, and the
seasonal periods of the underlying process are multiples of the
aggregation size. We derive the limit of the normalized spectral
density function of the aggregates, with increasing aggregation. The
limiting aggregate (seasonal) long-memory model may then be useful for
analyzing aggregate time-series data, which can be estimated by
maximizing the Whittle likelihood. We prove that the maximum Whittle
likelihood estimator (spectral maximum likelihood estimator) is
consistent and asymptotically normal, and study its finite-sample
properties through simulation. The efficacy of the proposed approach is
illustrated by a real-life internet traffic example.
\end{abstract}

%
\begin{keyword}
\kwd{asymptotic normality}
\kwd{consistency}
\kwd{seasonal auto-regressive fractionally integrated moving-average models}
\kwd{spectral density}
\kwd{spectral maximum likelihood estimator}
\kwd{Whittle likelihood}
\end{keyword}

\end{frontmatter}

\section{Introduction}
Data are often aggregated before analysis, for example,
1-minute data aggregated into half-hourly data or daily data aggregated
into monthly data.
Aggregation of data may be carried out for ease of interpretation
on a scale that is of interest, for example, policy makers and/or the public
are more interested in monthly unemployment rate than
daily unemployment rate.
On the other hand, data may be naturally aggregated, for example,
tree-ring data, which are often hard to disaggregate.
On a fine sampling scale, many time series are of long memory in the sense
that their spectral density functions admit a pole at the zero frequency.
A popular class of discrete time long memory processes are
autoregressive fractionally integrated moving
average (ARFIMA) models (see Granger and Joyeux~\cite{GraJoy80},
Hosking~\cite{Hos81}).
Man and Tiao~\cite{ManTia06} and Tsai and Chan~\cite{TsaCha05} showed
that temporal aggregation preserves the long-memory parameter
of the underlying ARFIMA process. Ohanissian, Russell and Tsay
\cite{OhaRusTsa08}
made use of this property in developing a test for long-memory.
Furthermore, as the extent of aggregation increases to infinity,
the limiting model retains the long-memory parameter of the original process,
whereas the short-memory components vanish.

In practice, the underlying process may admit seasonal long memory in
that its
spectral density function may have poles at certain non-zero frequencies.
Such data may be modeled as some
Seasonal Auto-Regressive Fractionally Integrated Moving-Average (SARFIMA)
process, see Section~\ref{sarfima}.
If the aggregation interval is much larger than the largest seasonal period,
aggregation will intuitively
merge the seasonal long-memory components with the regular long-memory
component and
eliminate the regular or seasonal short-memory components of the raw data.
For example, within the framework of ARIMA models,
Wei~\cite{Wei78} showed that aggregation
removes seasonality if the frequency of aggregation is larger than or
the same as the seasonal frequency.

On the other hand, if the aggregation
interval is large but is just some fraction
of the seasonal periods of the original data,
the aggregates may be expected to keep the seasonal
short- and long-memory pattern, albeit with different periods.
For many data, the latter scenario may be more relevant for analysis.
For example, aggregating 1-minute data into half-hourly data may remove the
short memory component on the minute scale but the daily or monthly
correlation pattern of the raw data may persist in the aggregates.

Here, our purposes are twofold. First, we study the
intermediate asymptotics of aggregating a SARFIMA process. In particular,
we derive the limiting (normalized)
spectral density function of an aggregated SARFIMA process
via the asymptotic framework where the seasonal
periods of the SARFIMA model are multiples of the aggregation interval and
the aggregation interval is large.
While the original time series is assumed to be a SARFIMA process, the limiting
result is robust to the exact form of
the short-memory and the regular long-memory components.
The limiting spectral density functions then
define a class of models suitable for
analyzing aggregate time series that may have regular or seasonal
long-memory and short-memory components.
Second, we derive the large-sample properties of the spectral maximum likelihood
estimator of the limiting aggregate SARFIMA model, obtained by
maximizing the Whittle likelihood.

The rest of the paper is organized as follows.
The SARFIMA model is reviewed in Section~\ref{sarfima}. In
Section~\ref{aggregate},
we derive the limiting spectral density function of an aggregate
SARFIMA process, under
the intermediate asymptotic framework.
Spectral maximum likelihood estimation of the limiting aggregate
SARFIMA model
and its large-sample properties are discussed in Section~\ref{est}. We
compare the empirical performance of the spectral maximum likelihood
estimator of the limiting model with that of the SARFIMA model by
Monte Carlo studies in Section~\ref{simulate}. The simulation results
suggest that fitting the limiting model to the aggregate data
generally reduces the bias in some long-memory parameters than simply
fitting a SARFIMA model.
We illustrate the use of the limiting
aggregate SARFIMA model and its possible gains in long-term forecasts
with a real application in Section~\ref{application}.
We conclude in Section~\ref{sconclude}. All proofs are collected in
the appendix
of Chan and Tsai~\cite{ChaTsa}.

\section{Seasonal autoregressive fractionally integrated moving
average models}
\label{sarfima}

We now briefly review the SARFIMA model which is widely useful in
scientific analysis; see Porter-Hudak~\cite{Por90}, Ray~\cite{Ray93},
Montanari, Rosso and Taqqu
\cite{MonRosTaq00}, Palma and Chan~\cite{PalCha05}, Bisognin and
Lopes~\cite{BisLop07} and
Lopes~\cite{Lop08}.
Let $\{Y_t,t=0,\pm1,\pm2,\ldots\}$ be a seasonal autoregressive
fractionally integrated moving average (SARFIMA) model with multiple
periods $s_1,\ldots,s_c$
%
\begin{equation}\label{eqn_sarfima}
\phi(B)(1-B)^{d}\prod_{i=1}^c\Phi_i(B^{s_i})(1-B^{s_i})^{D_i}Y_t
=\theta(B)\prod_{i=1}^c\Theta_i(B^{s_i})\varepsilon_t,
\end{equation}
where $d$ and $D_i, i=1,\ldots,c$, are real numbers,
$s_c>s_{c-1}>\cdots>s_1>1$ are integers, $\{\varepsilon_t\}$ is an
uncorrelated sequence of
random variables with zero mean and common, finite variance $\sigma
^2_\varepsilon>0$,
$\phi(z)=1-\phi_1z-\cdots-\phi_pz^p$, $\theta(z)=1+\theta
_1z+\cdots+\theta_qz^q$, and for $i=1,\ldots,c$,
$\Phi_i(z)=1-\Phi_{i,1}z-\cdots-\Phi_{i,P_i}z^{P_i}$, $\Theta
_i(z)=1+\Theta_{i,1}z+\cdots+\Theta_{i,Q_i}z^{Q_i}$,
$B$ is the backward shift operator, and $(1-B)^d$ is defined by the
binomial series expansion
\[
(1-B)^{d}= \sum_{k=0}^\infty\frac{\Gamma(k-d)}{\Gamma(k+1)\Gamma(-d)}B^k,
\]
where
$\Gamma(\cdot)$ is the gamma function.
Stationarity of $\{Y_t\}$ requires $D_i<1/2$ for all $i$ and
$d+\sum_{i=1}^c D_i< 1/2$, see Palma and Bondon~\cite{PalBon03}.
We assume that none of the roots of $\phi(\cdot)$ and
$\Phi_i(\cdot),$
$i=1,\ldots,c$,
match any roots of
$\theta(\cdot)$
and $\Theta_i(\cdot),$
$i=1,\ldots,c$.
Moreover, all roots of of the above polynomials are assumed to lie
outside the unit circle.
The conditions on the roots, the fractional orders $d$ and $D$'s ensure
that $\{Y_t\}$ is stationary and the model is identifiable.
It can be readily checked that
the spectral density of $\{Y_t\}$ equals, for $-\uppi< \omega\leq\uppi$,
\begin{eqnarray}\label{singularity}
h(\omega) 
&=&\frac{\sigma^2}{2\uppi}\biggl|\frac{\theta(\exp(\mathrm{i}\omega
))}{\phi(\exp(\mathrm{i}\omega))}\biggr|^2
\biggl|2\sin\biggl(\frac{\om}{2}\biggr)\biggr|^{-2\delta_0}
\prod_{j=1}^c\biggl|\frac{\Theta_j(\exp(\mathrm{i}s_j\omega
))}{\Phi_j(\exp(\mathrm{i}s_j\omega))}\biggr|^2\nonumber\\[-8pt]\\[-8pt]
&&{}\times\prod_{j=1}^c\prod_{k=1}^{\tau_j}
\bigl|\bigl(\exp(\mathrm{i}\nu_{jk})-\exp(\mathrm{i}\om)\bigr
)\bigl(\exp(-\mathrm{i}\nu_{jk})-\exp(\mathrm{i}\om)\bigr)\bigr
|^{-2\delta_{jk}},\nonumber
\end{eqnarray}
where $\delta_{0}=d+D_1+\cdots+D_c$; $\tau_j=[s_j/2]$, the greatest
integer $\leq s_j/2$;
$\nu_{jk}=2\pi k/s_j$, for $j=1,\ldots,c$, and $k=1,\ldots,\tau_j$;
$\delta_{jk}=D_j$, for $k=1,\ldots,\tau_j-1$,
$\delta_{j\tau_j}=D_j$ if $s_j=2\tau_j+1$, and $\delta_{j\tau
_j}=D_j/2$ if $s_j=2\tau_j$.
From (\ref{singularity}), we see that,
as $\omega\to0$, the spectral density $f(\omega)=\mathrm{O}(|\omega
|^{-2d-2D_1-\cdots-2D_c})$, whereas
for $j=1,\ldots,c$, $k=1,\ldots,\tau_j$, as $\om\to\nu_{jk}$,
$f(\om)= \mathrm{O}( |\omega-\nu_{jk}|^{-2D_j})$.
Given our interest in long-memory processes,
throughout this paper, the parameters $d$ and the $D_j$'s are restricted
by the inequality constraints: $0 \leq d+D_1+\cdots+D_c <1/2$, and
$0 \leq D_j < 1/2$, for $j=1,\ldots,c$.

\section{Aggregates of SARFIMA models}
\label{aggregate}

For non-stationary data, we assume that, after suitable regular and/or seasonal
differencing, the data become stationary and follow some stationary
SARFIMA model.
Specifically, let $r$ and $R_j$, $j=1,\ldots,c$, be non-negative integers and
$\{Y_t,t=0,\pm1,\pm2,\ldots\}$ a time series such that
$(1-B)^{r}(1-B^{s_1})^{R_1}\cdots(1-B^{s_c})^{R_c}Y_t$ is a
stationary SARFIMA model
defined by equation~(\ref{eqn_sarfima}). Therefore, $\{Y_t\}$
satisfies the difference equation
%
\begin{equation}\label{yt}
\phi(B)(1-B)^{r+d}\prod_{j=1}^{c}\Phi_j(B^{s_j})(1-B^{s_j})^{R_j+D_j}Y_t
=\theta(B)\prod_{j=1}^{c}\Theta_j(B^{s_j})\varepsilon_t,
\end{equation}
which is referred to as
the $\operatorname{SARFIMA}(p,r+d,q)\times(P_1,R_1+D_1,Q_1)_{s_{1}} \times\cdots
\times(P_c,R_c+D_c,Q_c)_{s_c}$ model.

Let $m \geq2$ be an integer and $X_{T}^{m}=\sum_{k=m(T-1)+1}^{mT}Y_{k}$
be the nonoverlapping $m$-temporal aggregates of $\{Y_t\}$. Let
$\nabla=1-B$ be the first difference operator, and $\nabla_s=1-B^s$
the lag-$s$ difference operator. Let
$R=(r,R_1,\ldots,R_c)$,
$\xi=(d;D_j,j=1,\ldots,c;\Phi_{i,j},i=1,\ldots,c,j=1,\ldots,P_i;\Theta
_{i,j},i=1,\ldots,c,j=1,\ldots,Q_i)$,
and assume $s_i=mz_i$, $i=1,\ldots,c$, where the $z_i$'s are positive integers.
Below we derive the spectral density
of the aggregates, and the limit of the normalized spectral densities with
increasing aggregation. The normalization that makes the spectral
densities integrate to 1 is necessary because, without normalization,
the variance of the aggregates
generally increases to infinity with increasing aggregation.

\begin{theorem}
\label{t1} Assume that $\{Y_t\}$ satisfies the difference equation
defined by \textup{(\ref{yt})}.
\begin{longlist}
\item[\textup{(a)}] For $ r \geq0$, $R_i\geq0$,
$i=1,\ldots,c$, and $m=2h+1$, the spectral density function of
$\{\nabla^{r}\nabla_{z_1}^{R_1}\cdots\nabla_{z_c}^{R_c}X_T^{m}\}$
is given by
%
\begin{eqnarray}\label{frdm}
f_{\xi,R,m}(\omega)
&=&\frac{1}{m}\biggl|2\sin\biggl(\frac{\omega}{2}\biggr)\biggr|^{2r+2}
\prod_{j=1}^c\biggl|2\sin\biggl(\frac{z_j\omega}{2}\biggr)\biggr|^{-2D_j}
\prod_{j=1}^c\biggl|\frac{\Theta_j(\exp(\mathrm{i}z_j\omega))}{\Phi_j(\exp
(\mathrm{i}z_j\omega))}\biggr|^2\no\\[-8pt]\\[-8pt]
&&{}\times\sum_{k=-h}^{h}\biggl|2\sin\biggl(\frac{\omega+2k\uppi}{2m}\biggr)\biggr|^{-2r-2d-2}
g\biggl(\frac{\omega+2k\uppi}{m}\biggr),\nonumber
\end{eqnarray}
where $g(\omega)=\sigma^2(2\uppi)^{-1}|\theta(\exp(\mathrm{i}\omega))|^2
|\phi(\exp(\mathrm{i}\omega))|^{-2}$ and $-\uppi<\om\leq\uppi$.

If $m=2h$, the spectral density is given by equation~\textup{(\ref{frdm})} with
the summation ranging from
$-h+1$ to $h$ for $-\uppi<\om\leq0$ and from $-h$ to $h-1$ for $0<\om
\leq\uppi$.
\item[\textup{(b)}] As $m\to\infty$, the normalized spectral density function of
$\{\nabla^{r}\nabla_{z_1}^{R_1}\cdots\nabla_{z_c}^{R_c} X_{T}^{m}\}
$ converges to
$f_{\xi,R}(\om)=K_{\xi,R}f^{*}_{\xi,R}(\om)$, where
%
\begin{eqnarray}\label{fr1}
f^{*}_{\xi,R}(\om)&=&\biggl|\sin\biggl(\frac{\om}{2}\biggr)\biggr|^{2r+2}\prod_{j=1}^c
\biggl|\sin\biggl(\frac{z_j\om}{2}\biggr)\biggr|^{-2D_j}\nonumber
\\[-8pt]\\[-8pt]
&&{}\times\prod_{j=1}^c\biggl|\frac{\Theta_j(\exp(\mathrm{i}z_j\omega))}{\Phi_j(\exp
(\mathrm{i}z_j\omega))}\biggr|^2
\sum_{k=-\infty}^{\infty}|\omega+2k\uppi|^{-2r-2d-2},\nonumber
\end{eqnarray}
where $K_{\xi,R}$ is the normalization constant ensuring that $\int
_{-\uppi}^{\uppi}f_{\xi,R}(\om)\,\mathrm{d}\om=1$.
\end{longlist}
\end{theorem}

\begin{remark}\label{r1}
The assumption that $s_i=mz_i$, for
$i=1,\ldots,c$, and $m\to\infty$ in Theorem~\ref{t1}(b)
should be interpreted as follows: the periodicities $s_i$'s are
multiples of the aggregation size $m$, and the
aggregation size is large. Consider two examples. Example~(1): hourly
data (that have a quarterly seasonality)
are aggregated into monthly data, so $s=2160$, $m=720$, and $z=3$, and
Example~(2): half-hourly data (that have a weekly seasonality)
are aggregated into daily data, so $s=336$, $m=48$, $z=7$.
\end{remark}

\begin{remark}\label{r2}
Note\vspace*{1pt} that $z_c>z_{c-1}>\cdots>z_1\geq1$.
For $j=1,\ldots,c$, $k=0,1,\ldots,[z_j/2]$, let $\om_{jk}=\nu_{j(mk)}=2\uppi
k/z_j$, then
both $m^{-2r-2d-1}f_{\xi,R,m}$ and $f_{\xi,R}$ are of order\break
$\mathrm{O}(|\omega|^{-2d-2D_1-\cdots-2D_c})$, for
$\om\to0$, and of order $\mathrm{O}(|\omega-\om_{jk}|^{-2D_j})$, for $\om
\to\om_{jk}$, $j=1,\ldots,c$, $k=1,\ldots,[z_j/2]$.
The above observations indicate that, if the periodicities $s_i$'s are
multiples of the aggregation size $m$,
then the aggregates and their limits preserve the long-memory and
seasonal long-memory parameters of the underlying SARFIMA
process, whereas the $z_i$'s become the periodicities of the aggregated
series.
\end{remark}

\begin{remark}\label{r3}
If $z_1=1$, the corresponding
seasonal long-memory component is confounded with the regular
long-memory component for the limiting aggregate process.
Hence, without loss of generality, we shall set $D_1=0$ if $z_1=1$ in
applications.
\end{remark}

\begin{remark}\label{r4}
If $r=0$, then the limiting model of the
aggregates of
$\{Y_{t}\}$ is simply a
$\operatorname{SARFIMA}(P_1, R_1+D_1,Q_1)_{z_{1}} \times\cdots\times
(P_c,R_c+D_c,Q_c)_{z_c}$ process with
fractional Gaussian noise as the driving noise process, where the
self-similarity parameter
(Hurst parameter) of the underlying fractional Gaussian process equals
$H=d+1/2$.
See Beran \cite{Ber94} for definition of the fractional Gaussian noise.
\end{remark}

\section{Spectral maximum likelihood estimator and its large sample properties}
\label{est}

We are interested in applying the long-memory limiting aggregate
process derived in Section~\ref{aggregate}
to data analysis. For this purpose, we assume (i) $0 \leq d+D_1+\cdots
+D_c <1/2$ and
(ii) $0 \leq D_j < 1/2$ for $j=1,\ldots,c$.
The limiting aggregate process is of long memory regularly or
seasonally if
either $0 < d+D_1+\cdots+D_c <1/2$ or $0 < D_j < 1/2$ for some
$j\in\{1,\ldots,c\}$. 
We also introduce the parameter $\sigma$ to account for the variance
of the data.
Furthermore, we assume $z_j$, $j=1,\ldots,c,$ are known.
Consider a time series $\{Y_i\}_{i=1-\delta}^N$, where $\delta$ is a
positive integer
to be defined below, such that, conditional on $\{Y_i\}_{i=1-\delta}^0$,
$\{\nabla^r\nabla_{z_1}^{R_1}\cdots\nabla_{z_c}^{R_c}Y_i\}_{i=1}^N$
is a stationary process with its spectral density defined by
%
\begin{equation}\label{spectralone}
f(\omega;\xi,R,\sigma^2)=\sigma^2f^{*}(\omega;\xi,R),
\end{equation}
where $f^{*}(\omega;\xi,R)$ is define in~(\ref{fr1}),
$\delta=\max_r+\sum_{i=1}^cz_i\cdot\max_{R_i}$; $\max_r$ and
$\max_{R_i}$, $i=1,\ldots,c$, are the largest possible values of
$r$ and
$R_i$, $i=1,\ldots,c$, respectively, which we will consider in simulation
studies and real data analysis in
Sections~\ref{simulate} and~\ref{application}. That is, the spectral
maximum likelihood estimators $\hat{r}$
and $\hat{R}_i$, $i=1,\ldots,c$, satisfy the conditions that $\hat{r}\in
\{0,\ldots,\max_r\}$
and $\hat{R}_i\in\{0,\ldots,\max_{R_{i}}\}$, for $i=1,\ldots,c$.

It can be easily checked that, conditional on $\{Y_i\}_{i=1-\delta}^0$,
the joint distributions of $\{\nabla^{r}\nabla_{z_1}^{R_1}\cdots
\nabla_{z_c}^{R_c} Y_{i}\}_{i=1}^N$
and $\{Y_i\}_{i=1}^N$ are the same.
Therefore, conditional on $\{Y_i\}_{i=1-\delta}^0$, the (negative)
log-likelihood
function of $\{Y_{i}\}$ can be approximated by the (negative) Whittle
log-likelihood function (see Hosoya~\cite{Hos96})
%
\begin{equation}\label{whittle}
-\tilde{l}(\xi,R,\sigma^2) = \sum_{j=1}^{T}
\biggl\{\log f(\omega_j;\xi,R,\sigma^2)+\frac{I_N(\omega_j;R)}{f(\omega
_j;\xi,R,\sigma^2)}\biggr\},
\end{equation}
where $\omega_j:=2\uppi j/N \in(0,\uppi)$ are the Fourier frequencies,
$T$ is the largest integer
$\leq(N-1)/2$, $I_N(\omega;R) = |\sum_{j=1}^N U_{j}(R)\exp(\mathrm{i}j\omega
)|^2/(2\uppi N)$,
and $U_i(R)=\nabla^r\nabla_{z_1}^{R_1}\cdots\nabla_{z_c}^{R_c}Y_i$,
$i=1,\ldots,N$.
In (\ref{whittle}), the computation of $f(\omega_j;\xi,R,\sigma^2)$
requires evaluation of an infinite sum.
Here, we adopt the method of Chambers~\cite{Cha96} to approximate $f(\omega
;\xi,R,\sigma^2)$ by
\begin{eqnarray*}\label{spectral2}
&&\tilde{f}(\omega;\xi,R,\sigma^2)\\
&&\quad=\sigma^2\biggl|\sin\biggl(\frac{\omega}{2}\biggr)\biggr|^{2r+2}
\prod_{j=1}^c \biggl|\sin\biggl(\frac{z_j\om}{2}\biggr)\biggr|^{-2D_j}
\prod_{j=1}^c\biggl|\frac{\Theta_j(\exp(\mathrm{i}z_j\omega))}{\Phi_j(\exp
(\mathrm{i}z_j\omega))}\biggr|^2
h(\omega;\xi,R),
\end{eqnarray*}
where $h(\omega;\xi,R) = \{2\uppi(2r+2d+1)\}^{-1}\{(2\uppi M-\omega
)^{-2r-2d-1}+(2\uppi M+\omega)^{-2r-2d-1}\}
+\sum_{k=-M}^M|\omega+2k\uppi|^{-2r-2d-2}$ for some large integer $M$.
By routine analysis, it can be shown that, under the conditions stated
in Theorem~\ref{asymptotic},
the approximation error of $h(\omega;R,\xi)$ to the infinite sum is
of order $\mathrm{O}(M^{-2r-2d-2})$.
Also, the approximation error of the first partial derivative with
respect to $d$ is of
order $\mathrm{O}(M^{-2r-2d-1-\epsilon})$, for any positive $\epsilon$ less
than~1.
These error rates guarantee that if the truncation parameter $M$
increases with
the sample size at a suitable rate, then the truncation has negligible
effects on the asymptotic distribution of the estimator, see
Theorem~\ref{asymptotic} below.
Replacing $f(\omega_j;\xi,R,\sigma^2)$ by $\tilde{f}(\omega_j;\xi
,R,\sigma^2)$ and letting
$\tilde{g}(\omega_j;\xi,R)=\tilde{f}(\omega_j;\xi,R,\sigma
^2)/\sigma^2$,
the (negative) Whittle log-likelihood function~(\ref{whittle}) now becomes
%
\begin{equation}\label{whittle1}
-\tilde{l}(\xi,R,\sigma^2) = \sum_{j=1}^{T}
\biggl\{\log\sigma^2 + \log\tilde{g}(\omega_j;\xi,R)
+\frac{I_N(\omega_j;R)}{\sigma^2\tilde{g}(\omega_j;\xi,R)}\biggr\}.
\end{equation}
Differentiating (\ref{whittle1}) with respect to $\sigma^2$ and
equating to zero gives
%
\begin{equation}\label{sigma_square}
\hat{\sigma}^2  =  \frac{1}{T}\sum_{j=1}^{T}\frac{I_N(\omega
_j;R)}{\tilde{g}(\omega_j;\xi,R)}.
\end{equation}

\noindent Substituting (\ref{sigma_square}) into (\ref{whittle1})
yields the objective function
%
\begin{equation}\label{obj}
-\tilde{l}(\xi,R)
 =  \sum_{j=1}^T \log\tilde{g}(\omega_j;\xi,R)
+T\log\Biggl(\sum_{j=1}^T\frac{I_N(\omega_j;R)}{\tilde{g}(\omega_j;\xi,R)}\Biggr)+C,
\end{equation}
where $C= T - T\log T$. The objective function is minimized with
respect to $\xi$ and $R$ to get the
spectral maximum likelihood estimators $\hat{\xi}$ and $\hat{R}$;
the estimator $\hat{\sigma}^2$ is then calculated by (\ref{sigma_square}).
Specifically, the spectral maximum likelihood estimators $\hat{\xi}$
and $\hat{R}$ are computed based on
equation~(\ref{obj}) using the following procedure
(Recall that $0 \leq D_j < 1/2$, for $j=1,\ldots,c$, and $0 \leq
d+D_1+\cdots+D_c <1/2$).
For each $\tilde{r}\in\{0,\ldots,\max_r\}$ and $R_i\in\{
0,\ldots,\max_{R_{i}}\}$, for $i=1,\ldots,c$,
we first find the local maximum likelihood estimator of $\xi$ in the
range that
$0 \leq D_j < 1/2$, $j=1,\ldots,c$, and $\tilde{r}\leq r+d+D_1+\cdots
+D_c < \tilde{r}+1/2$.
In our experiments, we let $\max_r=\max_{R_{1}}=\cdots
=\max_{R_{c}}=2$.
These local maximum likelihood estimators are then used to find the
global maximum likelihood estimator
of $\xi$ and $R=(r,R_1,\ldots,R_c)$.

For simplicity, let $\theta=(\xi,\sigma^2)$, $\hat{\theta}=(\hat
{\xi},\hat{\sigma}^2)$ and $\hat{R}$
be the spectral maximum likelihood estimator that minimizes
the (negative) Whittle log-likelihood function~(\ref{whittle1}).
Below, we derive the large-sample distribution of the spectral maximum
likelihood estimator.\looseness=1

\begin{theorem}
\label{asymptotic}
Let the data $Y=\{Y_{i}\}_{i=1}^N$ be such that
$\{\nabla^r\nabla_{z_1}^{R_1}\cdots\nabla_{z_c}^{R_c}Y_i\}_{i=1}^N$
is sampled from a stationary Gaussian seasonal long-memory process with the
spectral density given by \textup{(\ref{spectralone})}.
Let the spectral maximum likelihood estimator
$\hat{\theta}\in\Theta$, a compact parameter space, and the true
parameter $\theta_0$ be
in the interior of the parameter space.
Assume that each component of $R=(r,R_1,\ldots,R_c)$ is known to be
between $0$ and some integer $K$. Let $r_0$ and $d_0$ be the true
values of $r$ and $d$, and
the truncation parameter $M$ increase with the sample size so that
$M\to\infty$.
Then the spectral maximum likelihood estimator $\hat{R}$ and $\hat
{\theta}$ are consistent.
Moreover, if $\sqrt{N}M^{-2r_0-2d_0-1}\to0$ as $N\to\infty$, then
$\sqrt{N}(\hat{\theta}-\theta_0)$
converges in distribution to a normal random vector with mean 0 and
covariance matrix $\Gamma(\theta_0)^{-1}$ with\looseness=1
%
\begin{equation}\label{gamma1}
\Gamma(\theta) = \frac{1}{4\uppi}\int_{-\uppi}^{\uppi}(\updown\log
f(\omega;R,\theta))
(\updown\log f(\omega;R,\theta))'\,\mathrm{d}\omega,
\end{equation}\looseness=0
where $\updown$ denotes the derivative operator
with respect to $\theta$, and superscript $'$ denotes transpose.
\end{theorem}

\section{Empirical comparison between the limiting model and the
SARFIMA model}
\label{simulate}

Given aggregation is finite in practice,
fitting the limiting model~(\ref{spectralone}) to aggregate data may
result in bias, even though the bias vanishes with increasing aggregation.
On the other hand, ``to some extent, a discrete
time series model is conditional on the time scale,'' as remarked by a
referee. So, it is pertinent to compare the empirical performance of
the long-memory parameter estimators based on
the proposed limiting model with those based on
the SARFIMA model fitted to the aggregate data.
As aggregation carries a signature in the long-memory data structure as
spelt out in Theorem~\ref{t1},
fitting a SARFIMA model to aggregate data may
result in even larger bias on the long-memory parameters
than the limiting model.
Here, we report some simulation results for clarifying the
aforementioned issue.
Consider the aggregated time series $\{Y_i\}_{i=1-\delta}^N$ such that
$\{\nabla^r\nabla_{z}^{R}Y_i\}_{i=1}^N$
is a stationary process with its spectral density defined by
%
\begin{eqnarray}\label{frdm1}
f(\omega;r,d,D,\sigma^2)
&=&\sigma^2\biggl|\sin\biggl(\frac{\omega}{2}\biggr)\biggr|^{2r+2}
\biggl|\sin\biggl(\frac{z\omega}{2}\biggr)\biggr|^{-2D}\no\\[-8pt]\\[-8pt]
&&{}\times\sum_{k=-h}^{h-1}\biggl|(2m)\sin\biggl(\frac{\omega+2k\uppi}{2m}\biggr)\biggr|^{-2r-2d-2}
g\biggl(\frac{\omega+2k\uppi}{m}\biggr),\nonumber
\end{eqnarray}
where $g(\omega)=|\phi(\exp(\mathrm{i}\omega))|^{-2}$, $\phi(x)=1-\phi
_1x$, $\delta=\max_r+\max_{R}z$, and $0<\omega<\uppi$.
For $-\uppi<\omega\leq0$, $f(\omega;r,d,D,\sigma^2)=f(-\omega
;r,d,D,\sigma^2)$.
We consider $\sigma=2$, $z=10$, and $r=R=0$. The true values of $(d,
D)$ are (i) ($-0.1,0.3$) and (ii) ($0.2,0.25$), whereas
those of the other parameters are given in Table~\ref{tab21}.
The sample sizes considered are $N=512$ and $N=1024$.
We tried a range of $\operatorname{AR}(1)$ coefficient $\phi_1$: $-0.9, -0.5, 0.0,
0.5$, and $0.9$. The aggregation size $m$ are set to be $60$, $240$,
and $720$,
corresponding to the cases that minutely data are aggregated over
one hour, four hours, and half a day, respectively.
To each aggregated time series simulated from model~(\ref{frdm1}), we
fitted (i) the limiting aggregate model, and (ii) the SARFIMA model.
The averages and the standard deviations, as well as the asymptotic
standard errors, of 1000 replicates
of the estimators for (i)~$(d,D)=(-0.1,0.3)$ and (ii)
$(d,D)=(0.2,0.25)$ are summarized in Tables~\ref{tab21} and~\ref
{tab22}, respectively.
Note that the estimates of $r$ and $R$ equal \textit{zero} for all simulations.

\begin{sidewaystable*}
\tabcolsep=0pt
\tablewidth=\textwidth
\caption{Averages (standard deviations) of 1000 simulations of the spectral
maximum likelihood estimators of the parameters $d$, and $D$ by fitting
the limiting aggregate model~(\protect\ref{spectralone})
and the SARFIMA model~(\protect\ref{singularity}), respectively, to aggregate
data generated according to (\protect\ref{frdm1}), with $(d,D)=(-0.1,0.3)$.
The asymptotic standard errors for
the estimators of $(d,D,d+D)$ are $(0.03,0.03,0.04)$ and
$(0.02,0.02,0.03)$ for $N=512$ and $N=1024$, respectively.
Results under column heading ``A'' denote those from
the limiting aggregate model~(\protect\ref{spectralone}), whereas those under
``S'' are the counterparts from the SARFIMA model~(\protect\ref{singularity})}\label{tab21}
\begin{tabular*}{\textwidth}{@{\extracolsep{\fill}}llld{2.3}d{2.3}d{2.3}d{2.3}d{2.3}d{2.3}d{2.3}d{2.3}d{2.3}d{2.3}d{2.3}d{2.3}@{}} \hline
&&&\multicolumn{6}{l}{$N=512$}&\multicolumn{6}{l@{}}{$N=1024$}\\[-5pt]
&&&\multicolumn{6}{l}{\hrulefill}&\multicolumn{6}{l@{}}{\hrulefill}\\
 &\multirow{2}{*}{Para-} & \multicolumn{1}{l}{\multirow{2}{*}{True}} &\multicolumn{2}{l}{$m=60$}&\multicolumn{2}{l}{$m=240$}&\multicolumn{2}{l}{$m=720$}&\multicolumn{2}{l}{$m=60$}&\multicolumn{2}{l}{$m=240$}&\multicolumn{2}{l@{}}{$m=720$}
\\[-5pt]
&& &\multicolumn{2}{l}{\hrulefill}&\multicolumn{2}{l}{\hrulefill}&\multicolumn{2}{l}{\hrulefill}&\multicolumn{2}{l}{\hrulefill}&\multicolumn{2}{l}{\hrulefill}&\multicolumn{2}{l@{}}{\hrulefill} \\
\multicolumn{1}{l}{$\phi_1$}&meter & \multicolumn{1}{l}{value} & \multicolumn{1}{l}{A} & \multicolumn{1}{l}{S} & \multicolumn{1}{l}{A} & \multicolumn{1}{l}{S} & \multicolumn{1}{l}{A} &\multicolumn{1}{l}{S} & \multicolumn{1}{l}{A} &\multicolumn{1}{l}{S} & \multicolumn{1}{l}{A} &\multicolumn{1}{l}{S} & \multicolumn{1}{l}{A} & \multicolumn{1}{l@{}}{S}\\ \hline
$-$0.9&$d$ & $-$0.1 & -0.163 & -0.210 & -0.126 & -0.163 & -0.113& -0.147 & -0.162
& -0.206 & -0.125 & -0.160 & -0.112 & -0.144 \\
& & & (0.03) & (0.04) & (0.03) & (0.04) & (0.03)& (0.04) & (0.02) &
(0.03) & (0.02) & (0.03) & (0.02) & (0.03) \\
 &$D$ & \phantom{$-$}0.3 & 0.325 & 0.329 & 0.323 & 0.325 & 0.322 & 0.324 & 0.318 &
0.322 & 0.315 & 0.318 & 0.314 & 0.317 \\
& & & (0.04) & (0.04) & (0.04) & (0.04) & (0.04)& (0.04) &(0.03) &
(0.03) & (0.03) & (0.03) & (0.03) & (0.03) \\
&$d+D$ & \phantom{$-$}0.2 & 0.162 & 0.119 & 0.197 & 0.162 & 0.210 & 0.177 & 0.157 &
0.117 & 0.190 & 0.158 & 0.202 & 0.173 \\
& & & (0.05) & (0.05) & (0.05) & (0.05) & (0.05)& (0.05) & (0.03) &
(0.04) & (0.03) & (0.04) & (0.03) & (0.04) \\[3pt]
$-$0.5&$d$ & $-$0.1 & -0.112 & -0.146 &-0.105 & -0.138 & -0.103& -0.135 &-0.111 &
-0.143 &-0.105 & -0.135 &-0.103 & -0.133 \\
& & & (0.03) & (0.04) & (0.03) & (0.04) & (0.03)& (0.04) &(0.02) &
(0.03) &(0.02) & (0.03) &(0.02) & (0.03) \\
 &$D$ & \phantom{$-$}0.3 & 0.322 & 0.324 & 0.322 & 0.323 & 0.322 & 0.323 & 0.314 &
0.316 &0.314 & 0.316 &0.314 & 0.316 \\
& & & (0.04) & (0.04) & (0.04) & (0.04) & (0.04)& (0.04) & (0.03) &
(0.03) &(0.03) & (0.03) &(0.03) & (0.03) \\
&$d+D$ & \phantom{$-$}0.2 & 0.211 & 0.178 & 0.217 & 0.186 & 0.218 & 0.188 & 0.203 &
0.174 &0.209 & 0.180 & 0.211 & 0.183 \\
& & & (0.05) & (0.05) & (0.05) & (0.05) & (0.05)& (0.05) & (0.03) &
(0.04) &(0.03) & (0.04) & (0.03) & (0.04) \\[3pt]
\phantom{$-$}0&$d$ & $-$0.1 & -0.101 & -0.133 & -0.102 & -0.133 & -0.102& -0.134 &-0.100
& -0.131 &-0.101 & -0.131 &-0.101 & -0.131 \\
& & & (0.03) & (0.04) & (0.03) & (0.04) & (0.03)& (0.04) & (0.02) &
(0.03) &(0.02) & (0.03) &(0.02) & (0.03) \\
 &$D$ & \phantom{$-$}0.3 & 0.322 & 0.323 & 0.322 & 0.323 & 0.322 & 0.323 & 0.314 &
0.316 & 0.314 & 0.316 & 0.314 & 0.316 \\
& & & (0.04) & (0.04) & (0.04) & (0.04) & (0.04)& (0.04) & (0.03) &
(0.03) & (0.03) & (0.03) & (0.03) & (0.03) \\
&$d+D$ & \phantom{$-$}0.2 & 0.221 & 0.191 & 0.220 & 0.190 & 0.220 & 0.190 & 0.213 &
0.185 & 0.212 & 0.184 & 0.212 & 0.184 \\
& & & (0.05) & (0.05) & (0.05) & (0.05) & (0.05)& (0.05) & (0.03) &
(0.04) & (0.03) & (0.04) & (0.03) & (0.04) \\
\hline
\end{tabular*}
\end{sidewaystable*}

\setcounter{table}{0}
\begin{sidewaystable*}
\tabcolsep=0pt
\tablewidth=\textwidth
\caption{(Continued)}
\begin{tabular*}{\textwidth}{@{\extracolsep{\fill}}llld{2.3}d{2.3}d{2.3}d{2.3}d{2.3}d{2.3}d{2.3}d{2.3}d{2.3}d{2.3}d{2.3}d{2.3}@{}} \hline
&&&\multicolumn{6}{l}{$N=512$}&\multicolumn{6}{l@{}}{$N=1024$}\\[-5pt]
&&&\multicolumn{6}{l}{\hrulefill}&\multicolumn{6}{l@{}}{\hrulefill}\\
 &\multirow{2}{*}{Para-} & \multicolumn{1}{l}{\multirow{2}{*}{True}} &\multicolumn{2}{l}{$m=60$}&\multicolumn{2}{l}{$m=240$}&\multicolumn{2}{l}{$m=720$}&\multicolumn{2}{l}{$m=60$}&\multicolumn{2}{l}{$m=240$}&\multicolumn{2}{l@{}}{$m=720$}
\\[-5pt]
&& &\multicolumn{2}{l}{\hrulefill}&\multicolumn{2}{l}{\hrulefill}&\multicolumn{2}{l}{\hrulefill}&\multicolumn{2}{l}{\hrulefill}&\multicolumn{2}{l}{\hrulefill}&\multicolumn{2}{l@{}}{\hrulefill} \\
\multicolumn{1}{l}{$\phi_1$}&meter & \multicolumn{1}{l}{value} & \multicolumn{1}{l}{A} & \multicolumn{1}{l}{S} & \multicolumn{1}{l}{A} & \multicolumn{1}{l}{S} & \multicolumn{1}{l}{A} & \multicolumn{1}{l}{S} & \multicolumn{1}{l}{A} & \multicolumn{1}{l}{S} & \multicolumn{1}{l}{A} & \multicolumn{1}{l}{S} & \multicolumn{1}{l}{A} & \multicolumn{1}{l@{}}{S} \\ \hline
\phantom{$-$}0.5&$d$ & $-$0.1 & -0.092 & -0.121 & -0.099 & -0.130 &-0.101 & -0.132 &-0.091
& -0.119 & -0.098 & -0.128 & -0.100 & -0.130\\
& & & (0.03) & (0.04) & (0.03) & (0.04) &(0.03) & (0.04) &(0.02) &
(0.03) & (0.02) & (0.03) & (0.02) & (0.03) \\
 &$D$ & \phantom{$-$}0.3 & 0.321 & 0.323 & 0.322 & 0.323 & 0.322 & 0.323 &0.313 &
0.315 & 0.313 & 0.315 & 0.313 & 0.315 \\
& & & (0.04) & (0.04) & (0.04) & (0.04) & (0.04)& (0.04) &(0.03) &
(0.03) & (0.03) & (0.03) & (0.03) & (0.03) \\
&$d+D$ & \phantom{$-$}0.2 & 0.230 & 0.201 & 0.223 & 0.193 & 0.221 & 0.191 &0.222 &
0.195 & 0.215 & 0.187 & 0.213 & 0.186 \\
& & & (0.05) & (0.05) & (0.05) & (0.05) & (0.05)& (0.05) &(0.03) &
(0.04) & (0.03) & (0.04) & (0.03) & (0.04) \\[3pt]
\phantom{$-$}0.9&$d$ & $-$0.1 & -0.040 & -0.060 & -0.085 & -0.113 & -0.095& -0.126 &-0.041
& -0.059 & -0.085 & -0.111 &-0.095 & -0.124 \\
& & & (0.03) & (0.04) & (0.03) & (0.04) & (0.03)& (0.04) &(0.02) &
(0.03) & (0.02) & (0.03) &(0.02) & (0.03) \\
 &$D$ & \phantom{$-$}0.3 & 0.320 & 0.320 & 0.321 & 0.322 & 0.321 & 0.323 &0.311 &
0.311 & 0.313 & 0.314 &0.313 & 0.315 \\
& & & (0.04) & (0.04) & (0.04) & (0.04) & (0.04)& (0.04) &(0.03) &
(0.03) & (0.03) & (0.03) &(0.03) & (0.03) \\
&$d+D$ & \phantom{$-$}0.2 & 0.279 & 0.260 & 0.236 & 0.209 & 0.226 & 0.197 &0.270 &
0.252 & 0.228 & 0.203 &0.218 & 0.191 \\
& & & (0.05) & (0.05) &(0.05) & (0.05) & (0.05)& (0.05) &(0.03) &
(0.04) & (0.03) & (0.04) &(0.03) & (0.04) \\
\hline
\end{tabular*}
\end{sidewaystable*}

\begin{sidewaystable*}
\tabcolsep=0pt
\tablewidth=\textwidth
\caption{Averages (standard deviations) of 1000 simulations of the spectral
maximum likelihood estimators of the parameters $d$, and $D$ by fitting
the limiting aggregate model~(\protect\ref{spectralone})
and the SARFIMA model~(\protect\ref{singularity}), respectively, to aggregate
data generated according to (\protect\ref{frdm1}), with $(d,D)=(0.2,0.25)$.
The asymptotic standard errors for
the estimators of $(d,D,d+D)$ are $(0.03,0.03,0.04)$ and
$(0.02,0.02,0.03)$ for $N=512$ and $N=1024$, respectively.
Results under column heading ``A'' denote those from
the limiting aggregate model~(\protect\ref{spectralone}), whereas those under
``S'' are the counterparts from the SARFIMA model~(\protect\ref{singularity})}\label{tab22}
\begin{tabular*}{\textwidth}{@{\extracolsep{\fill}}lld{1.2}d{2.3}d{2.3}d{2.3}d{2.3}d{2.3}d{2.3}d{2.3}d{2.3}d{2.3}d{2.3}d{2.3}d{2.3}@{}}
\hline
&&&\multicolumn{6}{l}{$N=512$}&\multicolumn{6}{l@{}}{$N=1024$}\\[-5pt]
&&&\multicolumn{6}{l}{\hrulefill}&\multicolumn{6}{l@{}}{\hrulefill}\\
 &\multirow{2}{*}{Para-} & \multicolumn{1}{l}{\multirow{2}{*}{True}} &\multicolumn{2}{l}{$m=60$}&\multicolumn{2}{l}{$m=240$}&\multicolumn{2}{l}{$m=720$}&\multicolumn{2}{l}{$m=60$}&\multicolumn{2}{l}{$m=240$}&\multicolumn{2}{l@{}}{$m=720$}
\\[-5pt]
&& &\multicolumn{2}{l}{\hrulefill}&\multicolumn{2}{l}{\hrulefill}&\multicolumn{2}{l}{\hrulefill}&\multicolumn{2}{l}{\hrulefill}&\multicolumn{2}{l}{\hrulefill}&\multicolumn{2}{l@{}}{\hrulefill} \\
\multicolumn{1}{l}{$\phi_1$}&meter & \multicolumn{1}{l}{value} & \multicolumn{1}{l}{A} & \multicolumn{1}{l}{S} & \multicolumn{1}{l}{A} & \multicolumn{1}{l}{S} & \multicolumn{1}{l}{A} & \multicolumn{1}{l}{S} & \multicolumn{1}{l}{A} & \multicolumn{1}{l}{S} & \multicolumn{1}{l}{A} & \multicolumn{1}{l}{S} &\multicolumn{1}{l}{A} & \multicolumn{1}{l@{}}{S} \\ \hline
$-$0.9&$d$ & 0.2 & 0.186 & 0.223 & 0.197 & 0.231 & 0.198 & 0.232 & 0.189 &
0.225 & 0.198 & 0.234 & 0.200 & 0.235\\
& & & (0.03) & (0.03)& (0.03) & (0.03)& (0.03)& (0.03)& (0.02)& (0.02)&
(0.02)& (0.02)& (0.02) & (0.02)\\
 &$D$ & 0.25& 0.260 & 0.250 & 0.258 & 0.247 & 0.258 & 0.247 & 0.258 &
0.250 & 0.256 & 0.248 & 0.257 & 0.246 \\
& & & (0.04)& (0.03)& (0.04) & (0.03)& (0.04)& (0.03)& (0.03)& (0.02)&
(0.03)& (0.02)& (0.03) & (0.02)\\
&$d+D$ & 0.45& 0.446 & 0.474 & 0.455 & 0.479 & 0.456 & 0.479 & 0.446 &
0.476 & 0.454 & 0.482 & 0.457 & 0.481 \\
& & & (0.04)& (0.03)& (0.04) & (0.03)& (0.04)& (0.03)& (0.03)& (0.02)&
(0.03)& (0.02)& (0.03) & (0.02)\\[3pt]
$-$0.5&$d$ & 0.2 & 0.197 & 0.230 & 0.196 & 0.231 & 0.198 & 0.232 & 0.198 &
0.234 & 0.199 & 0.235 & 0.200 & 0.236\\
& & & (0.03)& (0.03)& (0.03) & (0.03)& (0.03)& (0.03)& (0.02)& (0.02)&
(0.02)& (0.02)& (0.02) & (0.02)\\
 &$D$ & 0.25& 0.258 & 0.247 & 0.258 & 0.247 & 0.258 & 0.248 & 0.257 &
0.248 & 0.256 & 0.246 & 0.257 & 0.246\\
& & & (0.04)& (0.03)& (0.04) & (0.03)& (0.04)& (0.03)& (0.03)& (0.02)&
(0.03)& (0.02)& (0.03) & (0.02) \\
&$d+D$ & 0.45& 0.455 & 0.478 & 0.454 & 0.479 & 0.456 & 0.480 & 0.454 &
0.481 & 0.455 & 0.482 & 0.457 & 0.483\\
& & & (0.04)& (0.03)& (0.04) & (0.03)& (0.04)& (0.03)& (0.03)&
(0.02)&(0.03) & (0.02)& (0.03) & (0.02)\\[3pt]
\phantom{$-$}0&$d$ & 0.2 & 0.198 & 0.232 & 0.199 & 0.231 & 0.199 & 0.231 & 0.200 &
0.236 & 0.199 & 0.235 & 0.200 & 0.236 \\
& & & (0.03)& (0.03)& (0.03) & (0.03)&(0.03) & (0.03)& (0.02)& (0.02)&
(0.02)& (0.02)& (0.02) & (0.02)\\
 &$D$ & 0.25& 0.258 & 0.247 & 0.258 & 0.246 &0.258 & 0.246 & 0.257 &
0.246 & 0.256 & 0.247 & 0.256 & 0.246 \\
& & & (0.04)& (0.03)& (0.04) & (0.03)&(0.04) & (0.03)&(0.03) & (0.02)&
(0.03)& (0.02)& (0.03) & (0.02)\\
&$d+D$ & 0.45& 0.456 & 0.479 & 0.457 & 0.477 & 0.456 & 0.478 & 0.457 &
0.483 & 0.456 & 0.482 & 0.456 & 0.482\\
& & & (0.04)& (0.03)& (0.04) & (0.03)& (0.04)& (0.03)& (0.03)& (0.02)&
(0.03)& (0.02)& (0.03) & (0.02)\\
\hline
\end{tabular*}
\end{sidewaystable*}

\setcounter{table}{1}
\begin{sidewaystable*}
\tabcolsep=0pt
\tablewidth=\textwidth
\caption{(Continued)}
\begin{tabular*}{\textwidth}{@{\extracolsep{\fill}}lld{1.2}d{2.3}d{2.3}d{2.3}d{2.3}d{2.3}d{2.3}d{2.3}d{2.3}d{2.3}d{2.3}d{2.3}d{2.3}@{}}
\hline
&&&\multicolumn{6}{l}{$N=512$}&\multicolumn{6}{l@{}}{$N=1024$}\\[-5pt]
&&&\multicolumn{6}{l}{\hrulefill}&\multicolumn{6}{l@{}}{\hrulefill}\\
 &\multirow{2}{*}{Para-} & \multicolumn{1}{l}{\multirow{2}{*}{True}} &\multicolumn{2}{l}{$m=60$}&\multicolumn{2}{l}{$m=240$}&\multicolumn{2}{l}{$m=720$}&\multicolumn{2}{l}{$m=60$}&\multicolumn{2}{l}{$m=240$}&\multicolumn{2}{l@{}}{$m=720$}
\\[-5pt]
&& &\multicolumn{2}{l}{\hrulefill}&\multicolumn{2}{l}{\hrulefill}&\multicolumn{2}{l}{\hrulefill}&\multicolumn{2}{l}{\hrulefill}&\multicolumn{2}{l}{\hrulefill}&\multicolumn{2}{l@{}}{\hrulefill} \\
\multicolumn{1}{l}{$\phi_1$}&meter & \multicolumn{1}{l}{value} & \multicolumn{1}{l}{A} & \multicolumn{1}{l}{S} & \multicolumn{1}{l}{A} & \multicolumn{1}{l}{S} & \multicolumn{1}{l}{A} & \multicolumn{1}{l}{S} & \multicolumn{1}{l}{A} & \multicolumn{1}{l}{S} & \multicolumn{1}{l}{A} & \multicolumn{1}{l}{S} & \multicolumn{1}{l}{A} & \multicolumn{1}{l@{}}{S}  \\ \hline
\phantom{$-$}0.5&$d$ & 0.2 & 0.201 & 0.233 & 0.199 & 0.233 &0.197 & 0.231 &0.202 & 0.238
&0.200 & 0.236 &0.200 & 0.235 \\
& & & (0.03) & (0.03)&(0.03) & (0.03)&(0.03) & (0.03)&(0.02) &
(0.02)&(0.02) & (0.02)&(0.02) & (0.02)\\
 &$D$ & 0.25& 0.257 & 0.245 & 0.258 & 0.246 &0.258 & 0.248 &0.256 &
0.245 &0.256 & 0.246 &0.256 & 0.247\\
& & & (0.04) & (0.03)& (0.04) & (0.03)&(0.04) & (0.03)&(0.03) &
(0.02)&(0.03) & (0.02)&(0.03) & (0.02) \\
&$d+D$ & 0.45& 0.459 & 0.478 & 0.457 & 0.479 &0.455 & 0.479 &0.458 &
0.483 &0.456 & 0.482 &0.456 & 0.483 \\
& & & (0.04)& (0.03)&(0.04) & (0.03)&(0.04) & (0.03)&(0.03) &
(0.02)&(0.03) & (0.02)&(0.03) & (0.02) \\[3pt]
\phantom{$-$}0.9&$d$ & 0.2 & 0.233 & 0.263 & 0.201 & 0.237 &0.200 & 0.233 &0.236 & 0.265
&0.206 & 0.241 & 0.201 & 0.237\\
& & & (0.03)& (0.03)&(0.03) & (0.03)&(0.03) & (0.03)&(0.02) &
(0.02)&(0.02) & (0.02)&(0.02) & (0.02)\\
 &$D$ & 0.25& 0.246 & 0.230 &0.254 & 0.244 &0.258 & 0.247 &0.247 &
0.231 &0.256 & 0.244 &0.256 & 0.246 \\
& & & (0.03)& (0.03)&(0.03) & (0.03)&(0.04) & (0.03)&(0.02) &
(0.02)&(0.02) & (0.02)&(0.03) & (0.02)\\
&$d+D$ & 0.45& 0.479 & 0.492 &0.455 & 0.481 &0.458 & 0.480 &0.489 & 0.495
&0.462 & 0.486 & 0.457 & 0.483\\
& & & (0.03)& (0.01)&(0.04) & (0.03)&(0.04) & (0.03)&(0.02) &
(0.01)&(0.03) & (0.02)&(0.03) & (0.02)\\
\hline
\end{tabular*}
\end{sidewaystable*}

From Tables~\ref{tab21} and~\ref{tab22},
it can be seen that the bias of the estimator of $d$ for
the limiting aggregate model is generally smaller in magnitude
than that of
the SARFIMA model,
except for $(d,D,\phi_1,m)=(-0.1,0.3,0.9,60)$,
and $(d,D,\phi_1,m)=(-0.1,0.3,0.9,240)$.
For $(d,D)=(0.2,0.25)$, the bias of the estimator of
$d+D$ for the limiting aggregate model is always smaller in magnitude
than that for the SARFIMA model.
For $(d,D)=(-0.1,0.3)$, the bias of $d+D$ for the limiting aggregate
model is smaller in magnitude than that for the SARFIMA model
if $\phi_1=-0.9$ or $\{(N,m,\phi_1)|N=1024,m=720,\phi_1\neq0.9\}$.
For $(d,D)=(-0.1,0.3)$, the bias of the estimator
of $D$ for the limiting aggregate model is always smaller than or equal
to that for the SARFIMA model, in magnitude.
For $(d,D)=(0.2,0.25)$, the bias of the estimator of
$D$ for the limiting aggregate model is always larger than that for the
SARFIMA model except for $(\phi_1,m)=(0.9,60)$
and $(\phi_1,m)=(0.9,240)$.
Overall, these limited simulation results suggest that the limiting model
leads to generally less biased estimates of $d$, and with comparable estimates
of $D$, than the SARFIMA model. Possible gains of long-term forecast accuracy
due to lesser bias in the estimator of $d$ based on the limiting model
will be further explored in the real application below.


\section{Application}
\label{application}

\begin{figure}
\begin{tabular}{@{}cc@{}}

\includegraphics{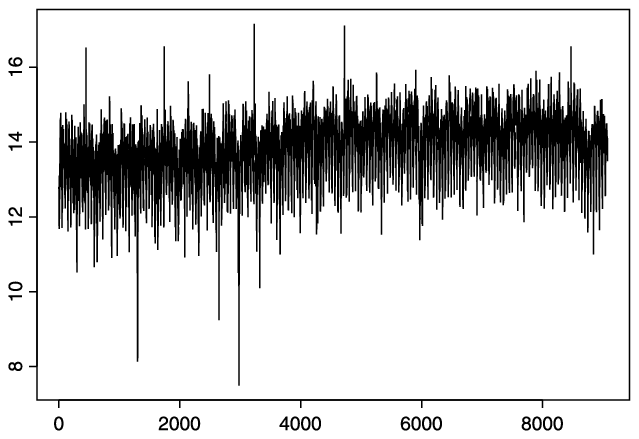}
&\includegraphics{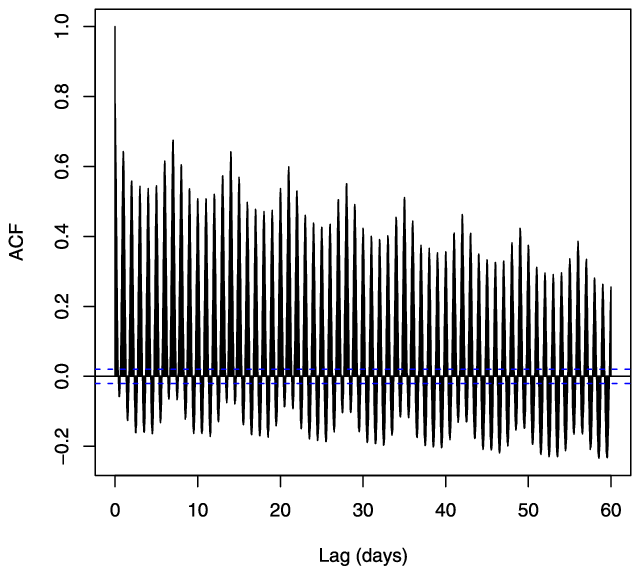}\\
(a)&(b)\\[6pt]
\multicolumn{2}{@{}c@{}}{
\includegraphics{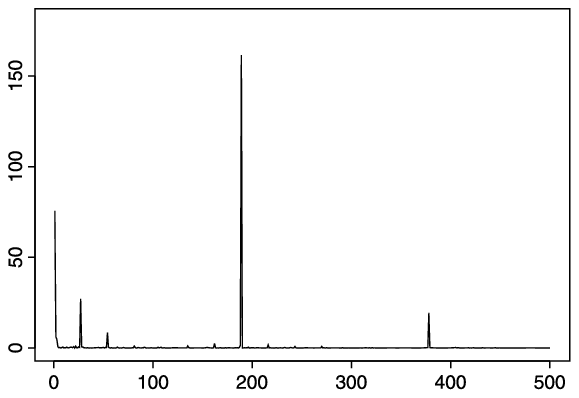}
}\\
\multicolumn{2}{@{}c@{}}{(c)}
\end{tabular}
\caption{Time series plot, sample ACF, and periodgram of the Log
transformed 30-minute (non-overlapping) aggregates. (a) Time series plot. (b) Sample ACF. (c) Periodogram.}\label{fig1}
\end{figure}

In this section, we report some analysis of
a time series of counts of http requests to a World Wide Web server at
the University of Saskatchewan, Canada,
between 1 June and 31 December in year 1995,
within the framework of the limiting aggregate seasonal long-memory model
and spectral maximum likelihood estimation.
The original data set consists of time stamps of 1-second
resolution, which can be downloaded from
\url{http://ita.ee.lbl.gov/html/contrib/Sask-HTTP.html}. Palma and
Chan~\cite{PalCha05} analyzed the
30-minute (non-overlapping) aggregates,
that is, each data point represents the total number of requests sent
to the Sakastchewan's server within a
30-minute interval. There are 9074 observations in total.
To make the data more Gaussian and to stabilize their variances, Palma
and Chan~\cite{PalCha05}
applied a logarithmic transformation to the aggregate data.
See Figure~\ref{fig1} for
the time series plot, the sample autocorrelation function, and the
periodogram of the transformed aggregate data.
Their fitted model is a $\operatorname{SARFIMA} (1,d,1) \times(0,D,0)_s$
model with $(\hat{d},\hat{D},\hat{\phi},\hat{\theta})=(0.076,0.148,0.917,0.583)$.
Although this model explains roughly two thirds of the total variance
of the data,
the residuals display significant autocorrelations at several lags,
in particular, at lags from 40 to 50 (Figure 6(a) of Palma and Chan~\cite{PalCha05}),
suggesting a lack of fit. Hsu and Tsai~\cite{HsuTsa09} also analyzed
the same data set, pointing out the presence of both daily and
weekly persistency in the data. Indeed, observe that there are two major
peaks in the periodogram: one at the origin and another at frequency
$\omega=2\uppi\times189/9074=0.1309$.
These features indicate a possible seasonal long-memory process with
$z=48$, that is, a daily pattern.
The third peak is at frequency $\omega=2\uppi\times27/9074=0.0187$,
indicating a possible weekly pattern.

Here, we reanalyze this dataset with the
limiting aggregate seasonal long-memory model defined by (\ref{spectralone})
with $c=3$, $r=R_1=R_2=R_3=0$, $z_1=1$, $z_2=48$ (corresponding to
daily effects),
and $z_3=48\times7=336$ (corresponding to weekly effects). Note that
$m=30\times60=1800$.
Our new approach may be justified as the 30-minute aggregation may well
fall within the intermediate asymptotic framework studied in
Section~\ref{aggregate}.
As discussed in Remark~\ref{r3} of Section~\ref{aggregate}, we assume $D_1=0$.
Specifically, if $\{Y_t\}$ is the observed time series,
the spectral density function of $\{Y_t\}$ can be written as
%
\begin{eqnarray}\label{fr_data}
\hspace*{-15pt}f_{\xi,\sigma^2}(\om)&=&\sigma^2\biggl|\sin\biggl(\frac{z_1\om}{2}\biggr)\biggr|^{2}
\biggl|\sin\biggl(\frac{z_2\om}{2}\biggr)\biggr|^{-2D_2}
\biggl|\sin\biggl(\frac{z_3\om}{2}\biggr)\biggr|^{-2D_3}
\biggl|\frac{\Theta_1(\exp(\mathrm{i}z_1\omega))}{\Phi_1(\exp(\mathrm{i}z_1\omega
))}\biggr|^2\no\\[-8pt]\\[-8pt]
\hspace*{-15pt}&&{}\times\biggl|\frac{\Theta_2(\exp(\mathrm{i}z_2\omega))}{\Phi_2(\exp
(\mathrm{i}z_2\omega))}\biggr|^2
\biggl|\frac{\Theta_3(\exp(\mathrm{i}z_3\omega))}{\Phi_3(\exp(\mathrm{i}z_3\omega))}\biggr|^2
\sum_{k=-\infty}^{\infty}|\omega+2k\uppi|^{-2d-2}.\nonumber
\end{eqnarray}

We have considered models of orders
$(P_1,Q_1,P_2,Q_2,P_3,Q_3)=(P_1,Q_1,0,0,0,0)$ with
$0\leq P_1\leq2$, and $0\leq Q_1\leq2$. The model with the smallest
AIC (Akaike information criterion)
is $(P_1,Q_1,P_2,Q_2,P_3,Q_3)=(2,2,0,0,0,0)$. Goodness of fit of this
model was studied in
Chan and Tsai~\cite{ChaTsa}.

The spectral maximum likelihood estimates of the parameters and the
$95\%$
bootstrap confidence intervals based on steps 1--4 of Section 6 of Chan
and Tsai~\cite{ChaTsa} are summarized in Table~\ref{tab_mle2200}.
The asymptotic standard deviations and the asymptotic $95\%$ confidence
intervals are also included in
Table~\ref{tab_mle2200}. It is clear that the bootstrap confidence
intervals of the parameters are comparable to their asymptotic
counterparts. The confidence intervals of
the parameters $d+D_2+D_3$, $D_2$ and $D_3$ indicate that the
long-memory pattern,
the daily seasonal long-memory pattern and the weekly seasonal
long-memory pattern are all significant.

\begin{table}
\caption{Spectral maximum likelihood estimates of
the parameters of the
model defined by equation~(\protect\ref{fr_data}), with
$(P_1,Q_1,P_2,Q_2,P_3,Q_3)=(2,2,0,0,0,0)$}\label{tab_mle2200}
\begin{tabular*}{\textwidth}{@{\extracolsep{\fill}}ld{2.4}lll@{}}
\hline
 &  & Bootstrap 95\% & Asymptotic&Asymptotic 95\%\\
& \multicolumn{1}{l}{Estimated} & confidence & standard & confidence \\
Parameter&\multicolumn{1}{l}{value} & interval & error &interval \\
\hline
$d$ & 0.2326 & \phantom{$-{}$}(0.1268, 0.2608) & 0.0436 & \phantom{$-{}$}(0.1471, 0.3181) \\
$D_2$ &0.1274 & \phantom{$-{}$}(0.1085, 0.1429) & 0.0083 & \phantom{$-{}$}(0.1111, 0.1437) \\
$D_3$ &0.1271 & \phantom{$-{}$}(0.1083, 0.1430) & 0.0083 & \phantom{$-{}$}(0.1108, 0.1434) \\
$d+D_2+D_3$ & 0.4871 & \phantom{$-{}$}(0.3821, 0.5000) & 0.0441 & \phantom{$-{}$}(0.4007, 0.5735) \\
$\phi_{1,1}$ &1.1277 & \phantom{$-{}$}(0.8916, 1.5089) & 0.1256 & \phantom{$-{}$}(0.8815, 1.3739) \\
$\phi_{1,2}$ &-0.2610 & ($-$0.5773, $-$0.0508) & 0.1009 & ($-$0.4588, $-$0.0632) \\
$\theta_{1,1}$ &-1.1788 & ($-$1.4586, $-$0.9315) & 0.0936 & ($-$1.3623, $-$0.9953) \\
$\theta_{1,2}$ &0.3593 & \phantom{$-{}$}(0.1237, 0.5755) & 0.0831 & \phantom{$-{}$}(0.1964, 0.5222)\\
$\sigma$ &0.3117 & \phantom{$-{}$}(0.3017, 0.3194) & 0.0051 & \phantom{$-{}$}(0.3017, 0.3217)\\
\hline
\end{tabular*}
\end{table}

\begin{figure}

\includegraphics{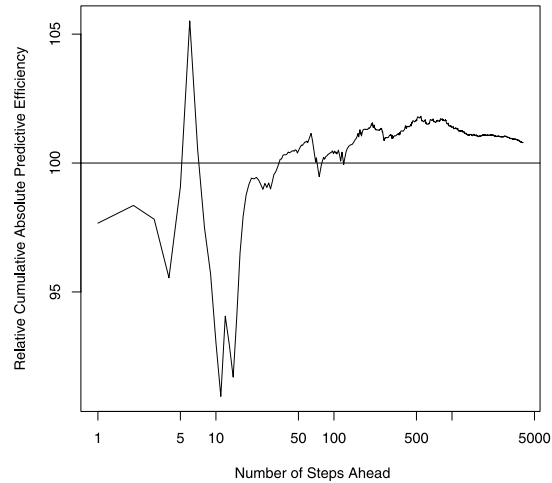}

\caption{Long-term forecast efficiency of the proposed model relative
to the SARFIMA model.}
\label{fig2}
\end{figure}

To assess the advantage of using the proposed aggregation model in
terms of forecasting,
as compared to a SARFIMA model, we divide the data roughly
into two halves, with the first 5074 data for fitting the proposed
model and
the SARFIMA model. We use the second
half of data for comparing their forecasting performance by
computing $h$-step ahead predictors, for $h=1,\ldots,4000$, and the
corresponding mean squared errors, via equations~(5.2.19) and (5.2.20)
of Brockwell and Davis~\cite{BroDav91}, respectively.
Specifically, the competing model we consider
is a $\operatorname{SARFIMA}(2,d,2)\times(0,D_1,0)_{s_{2}} \times(0,D_2,0)_{s_3}$
model, where
$s_2=48$, and $s_3=48\times7=336$. The estimates of the parameters are
$(\hat{d},\hat{D}_2,\hat{D}_3,\hat{d}+\hat{D}_2+\hat{D}_3,\hat
{\phi}_{1,1},\hat{\phi}_{1,2},
\hat{\theta}_{1,1},\hat{\theta}_{1,2},\hat{\sigma})
=(0.1996,0.1328,0.1280,0.4604,-0.038$, $0.8154,0.1099,-0.7342,0.1250)$.
Note that the estimators of $\hat{d}$ and $\hat{d}+\hat{D}_2+\hat{D}_3$
from the SARFIMA model are smaller than those from the proposed model, which
is similar to some finding reported in Section~\ref{simulate}.
Figure~\ref{fig2} displays
the ratios (in \%) of the cumulative first $h$ steps ahead mean
absolute forecast errors
of the SARFIMA model to their counterparts of
the limiting aggregate model, for $h=1,2,\ldots,4000$. These ratios
measure the long-term
forecast efficiency of the proposed model relative to the SARFIMA model.
As can be seen from the figure, the proposed model produces more
accurate long-term forecast
than the SARFIMA model once the forecast horizon is approximately
longer than 1 day.
We have also examined the rates the $h$-step ahead prediction variances
approach their
asymptotic value for the two models, and found that the
fitted proposed model
admits a slower convergence rate than the fitted SARFIMA model, which is
consistent with the longer memory (at the zero frequency) estimated by
the fitted proposed model than the SARFIMA model.

\section{Concluding remarks}
\label{sconclude}

We have derived the limiting structure of the temporal aggregates of a
(possibly non-stationary) SARFIMA model, with increasing aggregation.
We have also obtained some asymptotic properties of the spectral
maximum likelihood estimator of the limiting model,
including consistency and asymptotic normality. Monte Carlo experiments
show that the proposed method enjoys good empirical properties.
Moreover, estimator of $d$ under the proposed model appears to generally
have smaller bias than that from fitting a SARFIMA model to aggregate data.
The efficacy of our proposed methodology is illustrated with
an analysis of an internet traffic data.
Model diagnostic using a bootstrap procedure in the frequency domain, as
presented in Chan and Tsai~\cite{ChaTsa}, suggests a good fit.
Future research problems include extending the model to include
covariates and
developing other tools for model diagnostics.

\section*{Acknowledgements}
We are grateful to the referees for very helpful comments and thank
Academia Sinica, the National Science Council (NSC 95-2118-M-001-027), R.O.C.,
and the National Science Foundation (DMS-0405267, DMS-0934617) for
their support.

%

%

\printhistory

\end{document}